# Dirac matrices as elements of superalgebraic matrix algebra

## Monakhov V.V.

The paper considers a Clifford extension of the Grassmann algebra, in which operators are built from Grassmann variables and by the derivatives with respect to them. It is shown that a subalgebra which is isomorphic to the usual matrix algebra exists in this algebra, the Clifford extension of the Grassmann algebra is a generalization of the matrix algebra and contains superalgebraic operators expanding matrix algebra and produces supersymmetric transformations.

*УДК 530.145.63*

# МАТРИЦЫ ДИРАКА КАК ЭЛЕМЕНТЫ СУПЕРАЛГЕБРАИЧЕСКОЙ МАТРИЧНОЙ АЛГЕБРЫ

## В.В.Монахов

Санкт-Петербургский государственный университет, СПб, Россия

v.v.monahov@mail.ru

В работе рассмотрено клиффордово расширение алгебры Грассмана, в котором операторы построены из произведений грассмановых переменных и производных по ним. Показано, что в данной алгебре можно выделить подалгебру операторов, изоморфную обычной матричной алгебре, а сама алгебра является обобщением матричной алгебры, содержит суперальгебраические операторы, расширяющие матричную алгебру, и порождает преобразования суперсимметрии.

**Ключевые слова:** алгебра Грассмана, алгебра Клиффорда, обобщенная матричная алгебра, матрицы Дирака, суперпространство, суперсимметрия

Рассмотрим алгебру Грассмана $\Lambda_n$ [1]. В ней имеются подпространства ${}^k\Lambda$, где $k$ может принимать значения от 1 до $n$: ${}^1\Lambda$ – пространство мономов ранга 1 (линейное пространство с базисными элементами $\theta^\alpha$, где $\alpha = 1..n$), ${}^2\Lambda$ – мономов ранга 2 (с базисными элементами $\theta^\alpha \theta^\beta$, где $\alpha, \beta = 1..n, \alpha < \beta$), …, ${}^n\Lambda$ – мономов ранга $n$ [1]. Кроме того, добавим в алгебру единицу. Подпространство, натянутое на 1, будем называть $^0\Lambda$. Элементы этой алгебры

$$\Psi = \psi_0 + \psi_\alpha \theta^\alpha + \psi_{\alpha\beta}\theta^\alpha\theta^\beta + ... + \psi_{\alpha\beta...\eta}\theta^\alpha\theta^\beta...\theta^\eta \qquad (1)$$

будем рассматривать как векторы состояния, на которые могут действовать операторы.

На мономах и их линейных комбинациях вида (1) обычным образом [1] определим производные $\frac{\partial}{\partial\theta^\alpha}$ по антикоммутирующим переменным $\theta^\alpha$. Оператор умножения на $\theta^\alpha$ будем обозначать как $\hat{\theta}^\alpha$, оператор умножения на 1 будем обозначать как $\hat{1}$. Введем клиффордово расширение алгебры Грассмана [1,2] как алгебру $GM(\Lambda_n)$ [3] операторов с образующими $\hat{\theta}^\alpha$ и $\frac{\partial}{\partial\theta^\alpha}$. Эквивалентными будем считать операторы, результаты действия которых не отличаются при действии на произвольный вектор состояния (1).

Введем оператор $\hat{P}_n = \hat{\theta}^1\hat{\theta}^2...\hat{\theta}^n \frac{\partial}{\partial\theta^n}...\frac{\partial}{\partial\theta^2}\frac{\partial}{\partial\theta^1}$ – проектор на пространство $^n\Lambda$ мономов ранга $n$, затем оператор $\hat{P}_{n-1} = \hat{\theta}^{\alpha_1}\hat{\theta}^{\alpha_2}...\hat{\theta}^{\alpha_{n-1}}\frac{\partial}{\partial\theta^{\alpha_{n-1}}}...\frac{\partial}{\partial\theta^{\alpha_2}}\frac{\partial}{\partial\theta^{\alpha_1}}(\hat{1}-\hat{P}_n)$ - проектор на $^{n-1}\Lambda$, и так далее. В результате получим набор проекторов вплоть до $\hat{P}_1$ – проектора на $^1\Lambda$ (с базисными элементами $\theta^\alpha$) и $\hat{P}_0$ – проектора на $^0\Lambda$ (с базисным элементом 1).

Произвольный элемент $\xi \in {^1\Lambda}$ может быть представлен в виде

$$\xi = \xi_\alpha \theta^\alpha, \qquad (2)$$

а соответствующий ему элемент матричной алгебры будет записан в виде $\xi = \begin{pmatrix}\xi_1\\\xi_2\\...\end{pmatrix}$.

Сопоставим матрице $M$ произвольных линейных преобразований, переводящих элемент $\theta^\beta$ в элемент $\theta^\alpha$, элемент алгебры $GM(\Lambda_n)$ [3]

$$\hat{M} = m_\alpha{}^\beta \hat{\theta}^\alpha \frac{\partial}{\partial\theta^\beta} \hat{P}_1 \qquad (3)$$

Прямой проверкой легко убедиться, что произведение $\hat{A}\hat{B}$ операторов $\hat{A}$ и $\hat{B}$ вида

(2), соответствующих матрицам $\hat{A}$ и $\hat{B}$, соответствует произведению $AB$ матриц в матричной алгебре, и наоборот. То же относится к операциям умножения на число, сложения и комплексного сопряжения (если в соответствии с [1] считать, что комплексное сопряжение не меняет $\theta^\alpha$ и $\frac{\partial}{\partial \theta^\beta}$, и действует только на числовые коэффициенты).

Пусть $\lambda$ – число. Введем транспонированный к $\lambda\hat{\theta}^\alpha$ элемент как $(\lambda\hat{\theta}^\alpha)^T = \lambda\frac{\partial}{\partial\theta^\alpha}$, $(\lambda\frac{\partial}{\partial\theta^\alpha})^T = \lambda\hat{\theta}^\alpha$, и зададим по определению $(AB)^T = B^T A^T$. Легко видеть, что данные условия непротиворечивы. Легко заметить, что $(\hat{P}_n)^T = \hat{P}_n$ и $(\hat{1})^T = \hat{1}$. Поэтому $(\hat{P}_{n-1})^T = \hat{P}_{n-1}$, $(\hat{P}_{n-2})^T = \hat{P}_{n-2}$, …, $(\hat{P}_1)^T = \hat{P}_1$, $(\hat{P}_0)^T = \hat{P}_0$. Кроме того, $\hat{P}_1\hat{\theta}^\alpha = \hat{\theta}^\alpha\hat{P}_0$. Из чего следует выражение для транспонированной матрицы $M^T$:

$$(\hat{M})^T = (m^T)_\alpha{}^\beta \hat{P}_1\hat{\theta}^\alpha \frac{\partial}{\partial\theta^\beta} = (m^T)_\alpha{}^\beta \hat{\theta}^\alpha \hat{P}_0 \frac{\partial}{\partial\theta^\beta} = (m^T)_\alpha{}^\beta \hat{\theta}^\alpha \frac{\partial}{\partial\theta^\beta}\hat{P}_1.$$

То есть операция транспонирования элемента (2), соответствующего матрице, дает то же выражение, что и соответствие, сопоставляющее транспонированной матрице элемент алгебры $GM(\Lambda_n)$.

Таким образом, можно считать, что каждой матрице может быть взаимно однозначно сопоставлен элемент (3) алгебры, причем это соответствие сохраняется при умножении матриц на число, транспонировании матриц и их произведений, а также перемножении и сложении матриц. То есть алгебра элементов, заданных по формуле (3), и алгебра матриц изоморфны для всех матричных операций, включающих сложение, умножение и транспонирование матриц. Будем обозначать данную алгебру как $Mat(\Lambda_n)$ и называть ее супералгебраическим представлением матричной алгебры.

Обобщенная матричная алгебра $GM(\Lambda_n)$ шире обычной: помимо подалгебры $Mat(\Lambda_n)$ в ней существуют операторы обобщенных матриц [3]

$$\hat{M}(^{l}\Lambda, ^{k}\Lambda) = m_{\beta_1\ldots\beta_l}{}^{\alpha_k\ldots\alpha_1}\,\hat{\theta}^{\beta_1}\ldots\hat{\theta}^{\beta_l}\frac{\partial}{\partial\theta^{\alpha_k}}\ldots\frac{\partial}{\partial\theta^{\alpha_1}}\hat{P}_k, \qquad (4)$$

переводящие элементы из подпространства $^{k}\Lambda$ мономов ранга $k$ в подпространство $^{l}\Lambda$ мономов ранга $l$. Оператор-столбец является обобщенной матрицей $\hat{M}(^{1}\Lambda, ^{0}\Lambda) = m_{\alpha}\theta^{\alpha}\hat{P}_0$, переводящей $^{0}\Lambda$ в $^{1}\Lambda$, а оператор-строка – обобщенной матрицей $\hat{M}(^{0}\Lambda, ^{1}\Lambda) = m^{\beta}\frac{\partial}{\partial\theta^{\beta}}\hat{P}_1$, переводящей $^{1}\Lambda$ в $^{0}\Lambda$.

Рассмотрим теперь супералгебраическое представление матриц Дирака. Для начала рассмотрим ситуацию, когда число $n$ грассмановых переменных в столбце (2) равно порядку $n_D$ матриц Дирака. В этом случае в силу изоморфизма матричной алгебры и супералгебраического представления матричной алгебры матрицы Дирака имеют супералгебраическое представление в алгебре $Mat(\Lambda_n)$, а каждому столбцу-спинору можно сопоставить оператор-столбец обобщенной матричной алгебры. Поэтому все операторы алгебры $Mat(\Lambda_n)$ могут быть разложены по клиффордовым числам клиффордовой алгебры, в которой в качестве образующих выступают матрицы Дирака, а индексы $\alpha$ и $\beta$ в (3) являются спинорными индексами. При этом обобщенные матрицы (4) являются спин-тензорами.

В случае, когда число грассмановых переменных $n$ в столбце (2) превышает число $n_D$ элементов в столбце, соответствующем спинору Дирака в пространстве-времени, первые $n_D$ грассмановых переменных будем сопоставлять образующим подалгебры $Mat(\Lambda_{n_D})$, в которой матрицы Дирака в силу изоморфизма имеют супералгебраическое представление. Таким образом, $n_D$ грассмановых переменных оказываются образующими, порождающими матрицы Дирака и, следовательно, четные координаты суперпространства. При этом $n_{add} = n - n_D$ грассмановых переменных оказываются не связанными с матрицами Драка нечетными координатами суперпространства, будем называть их дополнительными грассмановыми образующими.

В алгебре $Mat(\Lambda_n)$ помимо операторов ее подалгебры $Mat(\Lambda_{n_D})$ имеются операторы, линейно преобразующие друг через друга последние $n_{add}$ грассмановых переменных, то есть не затрагивающие пространственные координаты, а также операторы, перемешивающие первые $n_D$ грассмановых переменных и последние $n_{add}$ грассмановых переменных. Будем в дальнейшем одним штрихом помечать индексы грассмановых образующих с номерами от 1 до $n_D$, а двумя штрихами – индексы дополнительных грассмановых образующих, с номерами от $n_D+1$ до $n_D+n_{add}$.

Если матрица Дирака $\Gamma^m$ в некотором представлении имеет матричные элементы $(\gamma^m)_\alpha{}^\beta$, в соответствии с ранее изложенным ей можно сопоставить обобщенную матрицу

$$\hat{\Gamma}^m = (\gamma^m)_{\alpha'}{}^{\beta'} \hat{\theta}^{\alpha'} \frac{\partial}{\partial \theta^{\beta'}} \hat{P}_1, \qquad (5)$$

где суммирование по индексам идет от 1 до $n_D$. При этом проектор $\hat{P}_1$ – это проектор из $Mat(\Lambda_n)$, а не из $Mat(\Lambda_{n_D})$, так как операторы $\hat{\Gamma}^m$ должны аннулировать все мономы, в которые входят множители $\theta^{\beta''}$.

Введем обозначения

$$\begin{aligned}\xi(n) &= \xi_{\alpha'}\theta^{\alpha'} + \xi_{\beta''}\theta^{\beta''} = \xi + \xi_{add}, \\ \xi'(n) &= \xi'_{\alpha'}\theta^{\alpha'} + \xi'_{\beta''}\theta^{\beta''} = \xi' + \xi'_{add}\end{aligned} \qquad (6)$$

Мы рассматриваем векторы состояния, поэтому в (6) пишем $\theta^{\alpha'}$ и $\theta^{\beta''}$ без шляпок. Рассмотрим теперь замену переменных: инфинитезимальные преобразования, перемешивающие первые $n_D$ грассмановых переменных и последние $n_{add}$ грассмановых переменных. В результате такого преобразования образующие основного базиса $\hat{\theta}^{\alpha'}$ преобразуются в

$$\begin{aligned}\tilde{\theta}^{\alpha'} &= \theta^{\alpha'} + \delta\theta^{\alpha'}, \\ \delta\theta^{\alpha'} &= \varepsilon^{\alpha'}_{\beta''}\theta^{\beta''}\end{aligned} \qquad (7)$$

где $\delta\theta^{\alpha'} = \varepsilon^{\alpha'}_{\beta''}\hat{\theta}^{\beta''}$ - изменения грассмановых образующих в результате преобразований, а $\varepsilon^{\alpha'}_{\beta''}$ - бесконечно малые числовые параметры. Здесь и далее мы пренебрегаем бесконечно малыми высших порядков.

Дополнительные грассмановы переменные $\tilde{\theta}^{\alpha''}$ преобразуются аналогичным образом:

$$\tilde{\theta}^{\beta''} = \theta^{\beta''} + \delta\theta^{\beta''},$$
$$\delta\theta^{\beta''} = \varepsilon^{\beta''}_{\alpha'}\theta^{\alpha'} \tag{8}$$

Поэтому

$$\tilde{\xi}(n) = \tilde{\xi}_{\alpha'}\tilde{\theta}^{\alpha'} + \tilde{\xi}_{\beta''}\tilde{\theta}^{\beta''} = (\tilde{\xi}_{\alpha'}\theta^{\alpha'} + \tilde{\xi}_{\alpha'}\varepsilon^{\alpha'}_{\beta''}\theta^{\beta''}) + (\tilde{\xi}_{\beta''}\theta^{\beta''} + \tilde{\xi}_{\beta''}\varepsilon^{\beta''}_{\alpha'}\theta^{\alpha'}) =$$
$$= (\tilde{\xi}_{\alpha'} + \varepsilon^{\beta''}_{\alpha'}\tilde{\xi}_{\beta''})\theta^{\alpha'} + (\tilde{\xi}_{\beta''} + \varepsilon^{\alpha'}_{\beta''}\tilde{\xi}_{\alpha'})\theta^{\beta''} \tag{9}$$

Формулы для $\xi'(n)$ аналогичны. Сами величины $\xi(n)$ и $\xi'(n)$ при замене переменных, естественно, не меняются, т.е. $\xi(n) = \tilde{\xi}(n)$ и $\xi'(n) = \tilde{\xi}'(n)$.

Поэтому $\xi_{\alpha'}\theta^{\alpha'} + \xi_{\beta''}\theta^{\beta''} = (\tilde{\xi}_{\alpha'} + \varepsilon^{\beta''}_{\alpha'}\tilde{\xi}_{\beta''})\theta^{\alpha'} + (\tilde{\xi}_{\beta''} + \varepsilon^{\alpha'}_{\beta''}\tilde{\xi}_{\alpha'})\theta^{\beta''}$.

Из чего следует $\xi_{\alpha'}\theta^{\alpha'} = (\tilde{\xi}_{\alpha'} + \varepsilon^{\beta''}_{\alpha'}\tilde{\xi}_{\beta''})\theta^{\alpha'}, \xi_{\beta''}\theta^{\beta''} = (\tilde{\xi}_{\beta''} + \varepsilon^{\alpha'}_{\beta''}\tilde{\xi}_{\alpha'})\theta^{\beta''}$.

Из линейной независимости величин $\theta^{\alpha}$ получаем, что коэффициенты $\xi_{\alpha}$ также преобразуются в бесконечно мало отличающиеся от первоначальных значения:

$$\tilde{\xi}_{\alpha'} = \xi_{\alpha'} + \delta\xi_{\alpha'}, \tilde{\xi}_{\beta''} = \xi_{\beta''} + \delta\xi_{\beta''}, \tag{10}$$

где, пренебрегая бесконечно малыми высшего порядка,

$$\delta\xi_{\alpha'} = -\varepsilon^{\beta''}_{\alpha'}\tilde{\xi}_{\beta''} = -\varepsilon^{\beta''}_{\alpha'}\xi_{\beta''}, \delta\xi_{\beta''} = -\varepsilon^{\alpha'}_{\beta''}\tilde{\xi}_{\alpha'} = -\varepsilon^{\alpha'}_{\beta''}\xi_{\alpha'}. \tag{11}$$

Полезно отметить, что в силу антикоммутации грассмановых переменных

$$\{\hat{\theta}^{\alpha'}, \delta\hat{\theta}^{\alpha'}\} = \{\hat{\theta}^{\alpha'}, \varepsilon^{\alpha'}_{\beta''}\hat{\theta}^{\beta''}\} = \varepsilon^{\alpha'}_{\beta''}\{\hat{\theta}^{\alpha'}, \hat{\theta}^{\beta''}\} = 0.$$

В соответствии с полученными ранее результатами зададим производные по величинам $\hat{\tilde{\theta}}^{\alpha'}$ как их транспонирование в обобщенной алгебре, и в результате получим

$$\frac{\partial}{\partial\tilde{\theta}^{\alpha'}} = \frac{\partial}{\partial\theta^{\alpha'}} + \varepsilon^{\alpha'}_{\beta''}\frac{\partial}{\partial\theta^{\beta''}}. \tag{12}$$

При этом с точностью до бесконечно малых высшего порядка для грассмановых переменных из основного базиса

$$\frac{\partial}{\partial \widetilde{\theta}^{\alpha'}}\widetilde{\theta}^{\gamma'} = (\frac{\partial}{\partial \theta^{\alpha'}} + \varepsilon_{\beta''}^{\alpha'}\frac{\partial}{\partial \theta^{\beta''}})(\hat{\theta}^{\gamma'} + \varepsilon_{\beta''}^{\gamma'}\hat{\theta}^{\beta''}) = \delta_{\alpha'}^{\gamma'} + \varepsilon_{\beta''}^{\gamma'}\delta_{\alpha'}^{\beta''} + \varepsilon_{\beta''}^{\alpha'}\delta_{\beta''}^{\gamma'} = \delta_{\alpha'}^{\gamma'},$$

что совпадает с соотношением $\frac{\partial}{\partial \theta^{\alpha'}}\theta^{\gamma'} = \delta_{\alpha'}^{\gamma'}$ для не преобразованных грассмановых переменных. Однако формулы для других производных по преобразованным переменным, вообще говоря, отличаются от формул для не преобразованных:

$$\frac{\partial}{\partial \widetilde{\theta}^{\alpha'}}\widetilde{\theta}^{\gamma''} = (\frac{\partial}{\partial \theta^{\alpha'}} + \varepsilon_{\beta''}^{\alpha'}\frac{\partial}{\partial \theta^{\beta''}})(\hat{\theta}^{\gamma''} + \varepsilon_{\beta'}^{\gamma''}\hat{\theta}^{\beta'}) = \delta_{\alpha'}^{\gamma''} + \varepsilon_{\beta'}^{\gamma''}\delta_{\alpha'}^{\beta'} + \varepsilon_{\beta''}^{\alpha'}\delta_{\beta''}^{\gamma''} = \varepsilon_{\alpha'}^{\gamma''} + \varepsilon_{\gamma''}^{\alpha'},$$

$$\frac{\partial}{\partial \widetilde{\theta}^{\gamma''}}\widetilde{\theta}^{\alpha'} = (\frac{\partial}{\partial \theta^{\gamma''}} + \varepsilon_{\beta'}^{\gamma''}\frac{\partial}{\partial \theta^{\beta'}})(\hat{\theta}^{\alpha'} + \varepsilon_{\beta''}^{\alpha'}\hat{\theta}^{\beta''}) = \delta_{\gamma''}^{\alpha'} + \varepsilon_{\beta''}^{\alpha'}\delta_{\gamma''}^{\beta''} + \varepsilon_{\beta'}^{\gamma''}\delta_{\beta'}^{\alpha'} = \varepsilon_{\gamma''}^{\alpha'} + \varepsilon_{\alpha'}^{\gamma''}$$

Потребуем, чтобы соотношения $\frac{\partial}{\partial \theta^{\alpha'}}\theta^{\gamma''} = 0$ и $\frac{\partial}{\partial \theta^{\gamma''}}\theta^{\alpha'} = 0$, означающие линейную независимость грассмановых переменных $\theta^{\alpha'}$ и $\theta^{\gamma''}$, выполнялись и для переменных $\widetilde{\theta}^{\alpha'}$ и $\widetilde{\theta}^{\gamma''}$, полученных в результате преобразования, т.е. чтобы было $\frac{\partial}{\partial \widetilde{\theta}^{\alpha'}}\widetilde{\theta}^{\gamma''} = 0$ и $\frac{\partial}{\partial \widetilde{\theta}^{\alpha''}}\widetilde{\theta}^{\gamma'} = 0$. Условия оказываются совместными, и из них следует

$$\varepsilon_{\alpha'}^{\gamma''} = -\varepsilon_{\gamma''}^{\alpha'}. \qquad (13)$$

Оператор $\hat{P}_1$ в (5) при данных преобразованиях не изменится, так как преобразование не меняет ранг мономов.

Рассмотрим пример подобного преобразования в простейшем случае, когда имеется три грассмановых образующих, две из которых, $\theta^1$ и $\theta^2$, соответствуют компонентам спинора трехмерного евклидова пространства, а третья, $\theta^3$ – дополнительная. Матрицы Дирака (гамма-матрицы) в данном пространстве – это хорошо известные матрицы Паули, которые, в соответствии с [4], обозначим как $H_m$, $m = 1,2,3$. Координаты $x_1, x_2, x_3$ вектора в этом пространстве можно получить как результат свертки двух спиноров $\xi = \begin{pmatrix} \xi_1 \\ \xi_2 \end{pmatrix}$ и $\xi' = \begin{pmatrix} \xi_1' \\ \xi_2' \end{pmatrix}$ [4]:

$$x_1 = (\xi')^T C H_1 \xi = \xi_1'\xi_1 - \xi_2'\xi_2,$$
$$x_2 = (\xi')^T C H_2 \xi = i(\xi_1'\xi_1 + \xi_2'\xi_2), \quad (14)$$
$$x_3 = (\xi')^T C H_3 \xi = -(\xi_1'\xi_2 + \xi_2'\xi_1),$$

где матрица $C = \begin{pmatrix} 0 & 1 \\ -1 & 0 \end{pmatrix}$.

В результате преобразования (7)-(13) получаются новые образующие $\tilde{\hat{\theta}}^{\alpha'}$ и производные по ним, и на основе замены в (5) первоначальных образующих на новые можно ввести новые матрицы Дирака и новые операторы-столбцы и операторы-строки. Мы получили новое трехмерное пространство, в котором эти матрицы являются клиффордовыми векторами, и в котором можно аналогично (14) ввести координаты

$$\tilde{x}_1 = (\tilde{\xi}')^T \tilde{C}\tilde{H}_1 \tilde{\xi},$$
$$\tilde{x}_2 = (\tilde{\xi}')^T \tilde{C}\tilde{H}_2 \tilde{\xi}, \quad (15)$$
$$\tilde{x}_3 = (\tilde{\xi}')^T \tilde{C}\tilde{H}_3 \tilde{\xi}$$

Тогда в соответствии с (14), (15) и (10), пренебрегая бесконечно малыми высшего порядка, получаем

$$\tilde{x}_1 = \tilde{\xi}_1'\tilde{\xi}_1 - \tilde{\xi}_2'\tilde{\xi}_2 = x_1 + \delta\xi_1'\xi_1 + \xi_1'\delta\xi_1 - \delta\xi_2'\xi_2 - \xi_2'\delta\xi_2 =$$
$$= x_1 + (\delta\xi')^T C H_1 \xi + (\xi')^T C H_1 \delta\xi,$$

и аналогичные формулы для $\tilde{x}_2$ и $\tilde{x}_3$. При этом из (11) имеем

$$\delta\xi_1 = -\varepsilon_1^3 \xi_3, \ \delta\xi_2 = -\varepsilon_2^3 \xi_3, \ \delta\xi_1' = -\varepsilon_1^3 \xi_3', \ \delta\xi_2' = -\varepsilon_2^3 \xi_3'.$$

В итоге получаем новое значение четных координат:

$$\tilde{x}_1 = x_1 + \delta x_1 = x_1 + (\delta\xi')^T C H_1 \xi + (\xi')^T C H_1 \delta\xi,$$
$$\tilde{x}_2 = x_2 + \delta x_2 = x_2 + (\delta\xi')^T C H_2 \xi + (\xi')^T C H_2 \delta\xi, \quad (16)$$
$$\tilde{x}_3 = x_3 + \delta x_3 = x_3 + (\delta\xi')^T C H_3 \xi + (\xi')^T C H_3 \delta\xi,$$

где

$$\delta\xi' = \delta\xi'_{\alpha'}\theta^{\alpha'} = -\varepsilon_{\alpha'}^3 \xi_3' \theta^{\alpha'},$$
$$\delta\xi = \delta\xi_{\alpha'}\theta^{\alpha'} = -\varepsilon_{\alpha'}^3 \xi_3 \theta^{\alpha'}, \quad (17)$$

а индекс $\alpha'$ меняется от 1 до 2.

В явном виде изменение четных координат:

$$\delta x_1 = \varepsilon_3^1(\xi_1'\xi_3 + \xi_3'\xi_1) - \varepsilon_3^2(\xi_2'\xi_3 + \xi_3'\xi_2),$$
$$\delta x_2 = i\varepsilon_3^1(\xi_1'\xi_3 + \xi_3'\xi_1) + i\varepsilon_3^2(\xi_2'\xi_3 + \xi_3'\xi_2), \quad (18)$$
$$\delta x_3 = -\varepsilon_3^1(\xi_2'\xi_3 + \xi_3'\xi_2) - \varepsilon_3^2(\xi_1'\xi_3 + \xi_3'\xi_1).$$

Из (7), (8), (13) и (12) находим изменение нечетных координат и новые значения производных по нечетным координатам:

$$\delta\theta^1 = \varepsilon_3^1\theta^3, \quad \delta\theta^2 = \varepsilon_3^2\theta^3, \quad \delta\theta^3 = \varepsilon_1^3\theta^1 + \varepsilon_2^3\theta^2, \quad \varepsilon_1^3 = -\varepsilon_3^1, \quad \varepsilon_2^3 = -\varepsilon_3^2, \quad (19)$$

$$\frac{\partial}{\partial\widetilde{\theta}^1} = \frac{\partial}{\partial\theta^1} + \varepsilon_3^1\frac{\partial}{\partial\theta^3}, \quad \frac{\partial}{\partial\widetilde{\theta}^2} = \frac{\partial}{\partial\theta^2} + \varepsilon_3^2\frac{\partial}{\partial\theta^3}.$$

Таким образом, любая функция $f(x_1, x_2, \theta^1, \theta^2)$, зависящая от четных переменных $x_1, x_2, x_3$ и нечетных $\theta^1, \theta^2$, может быть записана через новые координаты, полученные в результате преобразования (7)-(8), (11), (16)-(19):

$$f(\widetilde{x}_1, \widetilde{x}_2, \widetilde{x}_3, \widetilde{\theta}^1, \widetilde{\theta}^2) = f(x_1 + \delta x_1, x_2 + \delta x_2, x_3 + \delta x_3, \theta^1 + \delta\theta^1, \theta^2 + \delta\theta^2). \quad (20)$$

То есть рассмотренные инфинитезимальные преобразования грассмановых переменных $\theta^1$ и $\theta^2$ сопровождаются супертрансляцией (20), затрагивающей не только нечетные, но и четные координаты $x_1, x_2, x_3$. Данный эффект возникает из-за того, что матрицы (в том числе столбцы и строки) в (14)-(15) являются элементами обобщенной матричной алгебры и зависят от грассмановых переменных и производных по ним. При этом вектор состояния (1) можно считать суперполем – но, в отличие от имеющихся теорий, генераторы преобразования супертрансляции возникают как результат поворотов в пространстве нечетных координат, имеющих большую размерность, чем необходимо для формирования в суперагебраическом представлении матриц Дирака (гамма-матриц).

Проведенный анализ соответствует случаю N=1 суперсимметрии, когда дополнительных грассмановых переменных недостаточно для формирования еще одного комплекта матриц Дирака. Случаи N>1 требуют отдельного исследования.

Рассмотрение для случая четырех измерений должно быть аналогично проведенному для трех измерений, только вместо транспонированных спиноров в выражении, аналогичном (14), будут участвовать транспонированные комплексно сопряженные, а вместо $H_m$ будут использоваться $\Gamma_m$, $m = 0,1,2,3$.

Таким образом, матрицы Дирака являются супералгебраическими конструкциями – операторами обобщенной матричной алгебры. Что позволяет рассматривать матрицы Дирака как элементы обобщенной матричной алгебры, а использование операторов $M(^k\Lambda, ^l\Lambda)$ – как расширение спин-тензорного формализма на суперпространство. При количестве грассмановых переменных, превышающем число, необходимое для построения матриц Дирака, в рамках предложенного подхода естественным образом возникают преобразования суперсимметрии.

**Список литературы**